\documentclass[journal]{IEEEtran}

\usepackage{amsfonts}
\usepackage{epsfig}
\usepackage{subfig}
\usepackage{graphicx}
\usepackage{young}
\usepackage{enumerate}
\usepackage{url}

\newcommand{\beq}{\begin{equation}}
\newcommand{\eeq}{\end{equation}}
\newtheorem{thm}{Theorem}[section]

\newcommand{\field}[1]{\mathbb{#1}}

\begin{document}

\title{A phase-sensitive method for filtering on the sphere\thanks{This work was supported by an AcRF Tier 1 grant RG-033-09 to R. Kakarala.}}

\author{ Ramakrishna Kakarala and Philip Ogunbona \thanks{R. Kakarala is with the School of Computer Engineering, Nanyang Technological University, Singapore. Email: {\tt ramakrishna@ntu.edu.sg}; P. Ogunbona is with the School of Computer Science and Software Engineering, University of Wollongong, Australia. Emails: {\tt philipo@uow.edu.au}}}

\maketitle

\begin{abstract}
This paper examines filtering on a sphere, by first examining the roles of spherical harmonic magnitude and phase.  We show that phase is more important than magnitude in determining the structure of a spherical function.  We examine the properties of linear phase shifts in the spherical harmonic domain, which suggest a mechanism for constructing finite-impulse-response (FIR) filters.  We show that those filters have desirable properties, such as being associative, mapping spherical functions to spherical functions, allowing directional filtering, and being defined by relatively simple equations.  We provide examples of the filters for both spherical and manifold data.  
\end{abstract}

\begin{keywords} 
spherical harmonics, phase, FIR filtering
\end{keywords}


\thispagestyle{plain}
\markboth{For IEEE Trans. Signal Processing}{Phase-sensitive spherical filtering}
\section{Introduction}
While the importance of phase information in the ordinary Fourier transform on Euclidean domains $\field{R}^{n}$ is well-understood, it may be argued that a comparable understanding of phase does not exist when applied to the sphere $S^2$.  It is important to fill that gap because spherical data arise in numerous fields, including measurements of atmospheric pressure, cosmic microwave background, or surface reflectance.  Consequently, methods for filtering, spatio-spectral and spatio- scale analysis of spherical data using spherical harmonics have attracted the attention of the research community  \cite{DokmanicP10}\cite{ KhalidDSK12}\cite{McEwenHL08}.  However, existing research has not provided a detailed understanding of the importance of phase of spherical harmonics, either on its own or in relation to the harmonic magnitude.  It is the purpose of this paper to explore both the properties and the applications of phase information for spherical harmonics.   Our exploration of phase leads to a new method of filtering on the sphere which has desirable properties such as allowing for directional filtering, associative construction, and mapping spherical functions to spherical functions.  

\begin{figure}
 \centering
\subfloat[]{\includegraphics[width=0.6\linewidth]{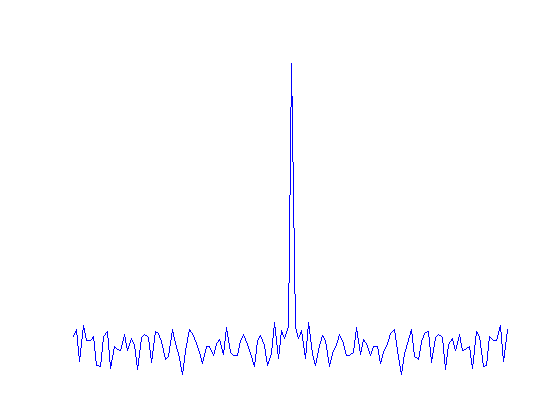}}
\subfloat[]{\includegraphics[width=0.6\linewidth]{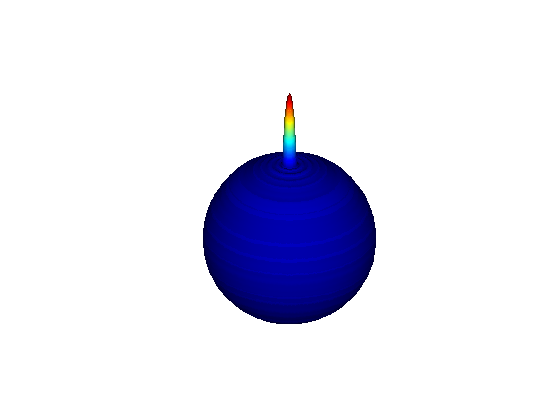}}
\caption{In part (a), the inverse Fourier transform of the magnitude-only spectrum from a white noise sample is shown.  The result is a positive definite function, which is an autocorrelation.  Part (b) shows the analogous result for the sphere, the inverse spherical harmonic expansion of the magnitude-only spectrum as discussed in this paper.}
\label{fig:posdef}
\end{figure}

We begin by recalling important properties of phase from the ordinary Fourier transform.  Let the Fourier transform of $f$ on the real line $\field{R}$ be denoted $F={\cal FT}\left\{f\right\}$.  The transform may be split into magnitude $|F|$ and phase $e^{j\phi}$ as follows: $F = |F|e^{j\phi}$.  If we set the phase to zero and invert the transform, i.e., by letting $G=|F|$ and $g ={\cal FT}^{-1}\left\{G\right\}$, then the resulting function $g$ must be a positive-definite function with its maximum value at the origin as illustrated in Figure~\ref{fig:posdef}(a).  In fact, it is well-known that $g$ is the autocorrelation of the function obtained from the inverse Fourier transform of the square-root magnitude, $\sqrt{|F|}$.   Recall that the autocorrelation is the integral
\beq
a_{f}(s) = \int_{\field{R}} f^{*}(x)f(x+s)dx.
\label{eq:acorr}
\eeq
Then, it is well known that $A_{f} = {\cal FT}\left\{a_{f}\right\} = FF^{*} = |F|^2$.  Hence, if its phase is set to zero, each function is turned into an auto-correlation.  Similarly, if we operate on photographic images, and swap the phase of image $f$ with that of image $g$, i.e., if we let $\widehat{f} = {\cal FT}^{-1}\left\{|F|\frac{G}{|G|}\right\}$, then it is well-known \cite{oppenheimwillsky} that $\widehat{f}$ resembles $g$ in appearance much more than $f$.  As we show below, similar demonstrations are possible with harmonic analysis on the sphere.  Figure~\ref{fig:posdef}(b) gives a preview by showing the result for $S^2$ that corresponds to removing phase entirely, using methods described later in this paper.   

Given the importance of phase in determining the appearance and structure of data, we might expect that phase information is crucial for applications such as viewpoint-invariant 3-D shape recognition in computer vision.  Yet, despite that, spherical harmonic invariants for shape recognition are in fact phase-blind magnitudes \cite{KazhdanFR03}\cite{ReisertB06}. A recent paper \cite{kakaralacvpr} showed that by using phase-sensitive invariants derived from the bispectrum of spherical harmonics, both discrimination and robustness to noise improve.    There has been little comment otherwise on the role that phase plays in determining spherical functions, which is one of the motivations for this paper.  

There is an extensive literature on filtering on the sphere. The mathematical aspects of that literature are discussed later, but here we summarize the contributions made to date and indicate how our approach is different.  A pioneering work by Driscoll \& Healy \cite{Driscoll_Healy_1994} describes convolution of two functions, $f$ and $h$, on the sphere by integrating the rotated version of $h$ with $f$.  The result is a function $g$ on $SO(3)$, which may be parameterized by the Euler angles denoted $\alpha$, $\beta$, $\gamma$ as $g(\alpha,\beta,\gamma)$.   For a fixed $\gamma$, the output $g_\gamma(\alpha,\beta)$ is also a function on $S^2$.   However, for the fixed $\gamma$, there is no harmonic space product expression that is  analogous to the familiar product $F\cdot G$ on $\field{R}$.   A product expression occurs, however, when $g$ is projected back to the sphere by integrating out one angle of rotation, which is equivalent to using rotationally-symmetric version  of the filter.     Those filtering methods are adopted in various papers \cite{Bulow04}\cite{DokmanicP10}\cite{McEwenHL08}\cite{McEwenHML07}\cite{WandeltGorski01}\cite{Wiauxeuclid}\cite{WiauxMNR}\cite{Wiaux06}\cite{YeoOG08}.  In this paper, we propose an alternative filtering method which has two main differences to the previous methods: it  uses a relatively simple harmonic-space product expression which works for directional filters, and it provides a means for constructing finite-impulse response (FIR) filters.  There are two key insights which lead to our results: the use of vector-matrix expressions that simplify construction of filters in the harmonic domain, and the analysis of phase in spherical harmonics.    

\section{Background material}  

The theory of $3$-D rotations, group representations, and spherical harmonics is well-described in various books \cite{chirikjian}\cite{hamermesh}\cite{hewittross}. We review key concepts and results briefly in this section. Every rotation is represented by three Euler angles, denoted $\alpha$, $\beta$, and $\gamma$, with $\alpha \in [0,2\pi]$ denoting rotation about the $z$ axis, $\beta \in [0,\pi]$ denoting rotation about the $y$ axis, and $\gamma \in [0,2\pi]$ also denoting rotation about the $z$ axis.  Two of the angles, $\alpha$ and $\beta$, parameterize the sphere $S^2$, with $0\leq\alpha\leq 2\pi$ representing longitude (angle measured counter-clockwise from the positive $x$-axis in the $x$-$y$ plane), $0\leq\beta\leq\pi$
the colatitude (angle with respect to positive $z$-axis).  Each unit vector $u$ lying on $S^{2}$ may be described parametrically as 
\beq
u = \left[ \cos(\alpha)\sin(\beta), \sin(\alpha)\sin(\beta),\cos(\beta)\right]^{\top}.
\label{eq:unitvec}
\eeq
The Laplace spherical harmonics form an orthonormal basis for functions on $S^{2}$.  For each non-negative integer $\ell$, and integer $-\ell \leq m \leq \ell$ they have the functional form
\beq
Y_{\ell}^{m}(u) = Y_{\ell}^{m}(\beta,\alpha) = c_{\ell}{^m} P_{\ell}^m(\cos\beta)e^{-jm\alpha}, \quad -\ell\leq m\leq \ell.
\label{eq:sph}
\eeq
Here $P_{\ell}^m$ are the {\em associated Legendre functions}, which are real-valued, and the normalization constants are
\beq
c_{\ell}^{m} = \sqrt{\frac{2\ell+1}{4\pi}\frac{(\ell-m)!}{(\ell+m)!}}.
\label{eq:alpha}
\eeq
For convenience, we denote the integral on the sphere as follows
\beq
\int_{S^{2}} f(u) du = \int_{0}^{2\pi} \int_{0}^{\pi} f(\beta,\alpha) \sin (\beta) d\beta\, d\alpha.
\label{eq:sphere}
\eeq
With the above definitions, we have orthonormality
\beq
\int_{S^{2}} Y_{\ell}^{m}(u) Y_{p}^{n}(u)^{*} du = \delta_{\ell p}\delta_{m n}.
\label{eq:normal}
\eeq
Defining the $\ell$, $m$-th spherical harmonic coefficient as 
\beq
F_{\ell}^{m} = \int_{S^{2}} f(u) Y_{\ell}^{m}(u)^{*} du,
\label{eq:Fln}
\eeq
every square-integrable function $f$ on the sphere may be expanded as follows:
\beq
f(u)  = \sum_{\ell = 0}^{\infty}\sum_{m=-\ell}^{\ell} F_{\ell}^{m} Y_{\ell}^{m}(u)
\label{eq:spharmexp}
\eeq

It proves very convenient to write (\ref{eq:spharmexp}) as a series of vector inner products.  The vector notation not only reduces the clutter of indices but also provides linear-algebraic insights which led to the main results in this paper.   First, let us combine all coefficients for a given $\ell$ into a $1\times (2\ell + 1)$ row vector
\beq
F_{\ell} = \left[F_\ell^{-\ell},\ldots,F_{\ell}^{0},\ldots,F_{\ell}^{\ell}\right].
\label{eq:fourvect}
\eeq
We refer to $F_\ell$ below as the $\ell$-th Fourier coefficient of $f$.   Second, let $Y_{\ell}$ denote the $(2\ell + 1)\times 1$ {\em column} vector containing the spherical harmonics at the $\ell$-th frequency
\beq
Y_{\ell}(u) = \left[Y_{\ell}^{-\ell}(u),\ldots,Y_{\ell}^{0}(u), \ldots,Y_{\ell}^{\ell}(u) \right]^{\top}.  
\label{eq:yvec}
\eeq
Using these vectors, we write (\ref{eq:spharmexp}) as 
\beq
f(u)  = \sum_{\ell = 0}^{\infty} F_{\ell} Y_{\ell}(u)
\label{eq:vecnotspharmexp}
\eeq
This relatively-simple equation is, as we see below, the key to understanding the role of phase and the implementation of filtering.  Suppressing indices, we may ''visualize'' the equation as a sum of inner products of progressively longer odd-length vectors:
\beq
f(u) = F Y + \left[ \begin{array}{c} F\,F\,F\end{array} \right] \left[\begin{array}{c}Y \\ Y \\ Y\end{array}\right] +\left[ \begin{array}{c} F\,F\,F\,F\,F\end{array} \right] \left[\begin{array}{c}Y \\ Y \\ Y\\ Y\\ Y\end{array}\right] \ldots
\eeq
Here the symbols $F$ and $Y$ inside the vectors denote appropriate elements of the $F_\ell$ and $Y_{\ell}$ vectors with their indices suppressed.  
The same expansion (\ref{eq:vecnotspharmexp}), using only a finite bandwidth $L$ is 
\beq
f_{L}(u)  = \sum_{\ell = 0}^{L-1} F_{\ell} Y_{\ell}(u). 
\label{eq:finbandspharmexp}
\eeq
%
%

\subsection{Magnitude and phase}

It is interesting at this point to consider what might be the ``phase'' of spherical harmonic coefficients.  We may, for example, split each scalar
coefficient $F_{\ell}^{m}$ into magnitude $|F_{\ell}^{m}|$ and ``phase'' $e^{j\phi(\ell,m)}=F_{\ell}^{m}/|F_{\ell}^{m}|$.   However, that phase
does not behave analogously to the phase spectrum of the ordinary Fourier transform.  When a function $f$ on $\field{R}$ is translated 
by $x\mapsto x+t$, then its phase spectrum transforms as $\phi\mapsto \phi+2\pi \nu t$, where $\nu$ is the frequency.  There is no similar behavior for the scalar ``phase'' $\phi(\ell,m)$ of spherical harmonic coefficients just defined.

To determine a more appropriate interpretation of phase for the spherical harmonics, we examine their rotation property.  Every $3$-D rotation is
represented by $3\times 3$ orthogonal matrix $R$ with determinant $+1$.  If, for some $R$, we have $g(u) = f(Ru)$ for all $u\in S^2$, so that $g$ is a rotated version of $f$, then the Fourier coefficients of $g$ are related to those of $f$ as follows:
\beq
G_{\ell} = F_{\ell}D_{\ell}(R).
\label{eq:rotS2}
\eeq
Here $D_{\ell}(R)$ is a $(2\ell +1)$-dimensional unitary matrix, known as the Wigner $D$-matrix, with the homomorphic property 
$D_{\ell}(RS) = D_{\ell}(R)D_\ell(S)$ for every pair of rotations $R$, $S$.  The elements of the $D$-matrices are separable in
the Euler angles:
\beq
D_{\ell}^{mn}(\alpha,\beta,\gamma) = e^{-jm\alpha}d_{\ell}^{mn}(\beta)e^{-jn\gamma}, \quad -\ell \leq m,n \leq \ell.
\label{eq:wignerdfunct}
\eeq
Here $d_{\ell}^{mn}$ is the ``little'' Wigner $d$-function, which is real-valued for the $z$-$y$-$z$ Euler angles as defined earlier.  In particular, the spherical harmonics are, up to a scale factor, the elements of the {\em middle} column of the corresponding $D$ matrix:
\beq
Y_{l}^{m}(\beta,\alpha) = c_{\ell}^{0} D_{l}^{m0}(\alpha,\beta,0), \quad -\ell \leq m \leq \ell,
\label{eq:sphd}
\eeq
where, from (\ref{eq:alpha}), we have the constants
\beq
c_{\ell}^{0} = \sqrt{\frac{2\ell+1}{4\pi}}.
\label{eq:cellzero}
\eeq

We now obtain an appropriate definition of phase. Each coefficient vector $F_{\ell}$ may be separated into magnitude $\|F_{\ell}\|$ and ``phase'' determined 
by the unit-length vector $U_{\ell} = F_{\ell}/\|F_{\ell}\|$, so that $F_{\ell}=\|F_{\ell}\|U_{\ell}$.  Under a rotation of the underlying
function on the sphere, each coefficient $F_\ell$ transforms by (\ref{eq:rotS2}).  Since the Wigner matrices are unitary, we see that the magnitude
$\|F_{\ell}\|$ will not change, but the phase vector transforms as $U_{\ell}\mapsto U_{\ell}D_{\ell}(R)$.  That behavior matches the linear shift 
property of the ordinary Fourier transform under translation.  


\section{Properties of magnitude and phase}
\label{sec:magphase}

Any attempt to define magnitude and phase is made clearer with an exploration of properties.  We proceed by looking to the familiar properties on the real line $\field{R}$.  There, the Fourier transform of the Dirac delta function is constant for all frequencies:
\beq
F(\omega) = \int_{-\infty}^{\infty} \delta (x) e^{-j\omega x} dx = 1.
\eeq
When the delta function moves to any location $x_0$, only the phase of $F$ changes: $F(\omega) = e^{-j\omega x_0}$.

Unlike $\field{R}$, it is not obvious how to define a delta function on $S^{2}$, since the sphere does not have an ``origin'' with unique properties.  We take our inspiration from the equation
\beq
\delta(x-y) = \int_{-\infty}^{\infty} e^{-j\omega x}e^{j\omega y} \frac{d\omega}{2\pi}.
\label{eq:deltar}
\eeq
Following \cite[pp 594-595]{riley}, we define for any two points $u$, $w$ on $S^{2}$, the function (with $\dagger$ denoting conjugate-transpose) 
\beq
\delta(u-w) = \sum_{\ell=0}^{\infty} Y_{\ell}(u)^{\dagger} Y_{\ell}(w)
\label{eq:deltadef}
\eeq
It is easy to show, for any function $f$, that 
\beq
\int_{S^{2}} \delta(u-w) f(u) du  = f(w).
\label{eq:deltaprop}
\eeq

\begin{figure}
\includegraphics[width=\linewidth]{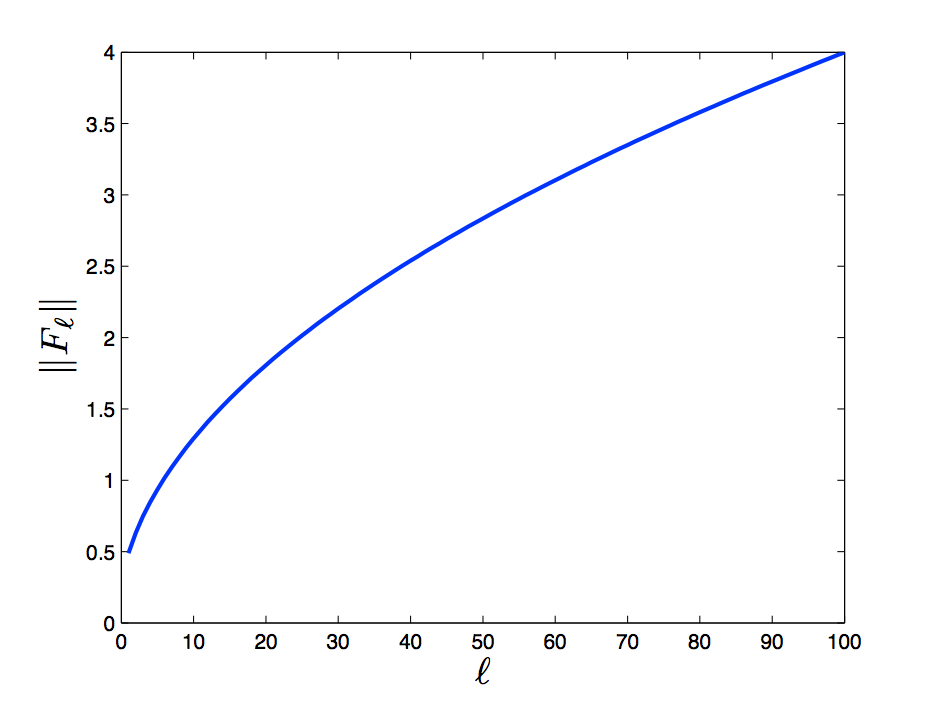}
\caption{The magnitude spectrum of the $\delta$ function on the sphere $S^{2}$. It increases with frequency $\ell$ unlike the counterpart on the real line $\field{R}$.}
\label{fig:deltaspec}
\end{figure}

Suppose now that we have a delta function $\delta(u-n)$ located at the North pole $n=[0,0,1]^{\top}$ of the sphere.  We take the angular coordinates of the north pole to be $(\beta,\alpha) = (0,0)$.  The spherical harmonic coefficients are 
\beq
\int_{S^{2}} \delta(u-n) Y_{\ell}^{m}(u)^{*} du = Y_{\ell}^{m}(0,0)^{*} = c_{\ell}^{m} P_{\ell}^{m}(1).
\eeq
Since $P_{\ell}^{m}(1) = 1$ for $m=0$, and zero for other $m$, we have that the Fourier coefficients of the delta function are
\beq
F_{\ell}^{m}  = \left\{ \begin{array}{rl}
 c_{\ell}^{0} &\mbox{ if $m=0$} \\
  0 &\mbox{ otherwise}
       \end{array} \right.
       \label{eq:impfourier}
\eeq
Consequently, the vectors $F_{\ell}$ are all zero except in the middle entry, which is equal to $c_\ell$.  The magnitude spectrum of the delta function is independent of its position on the sphere since it is rotation invariant.  Therefore, for a delta function located anywhere on the sphere, we have 
\beq
\|F_{\ell} \| = c_{\ell}^{0} = \sqrt{\frac{2\ell+1}{4\pi}}.
\label{eq:frac}
\eeq
This result has been published before \cite[eq 64]{ramamoorthi}.  Therefore, unlike on the real-line,  the magnitude spectrum of a delta function on $S^{2}$ is not constant, but rather increases with frequency $\ell$ in the manner illustrated by Figure~\ref{fig:deltaspec}.

\begin{figure}
\begin{tabular}{c}
\epsfig{file=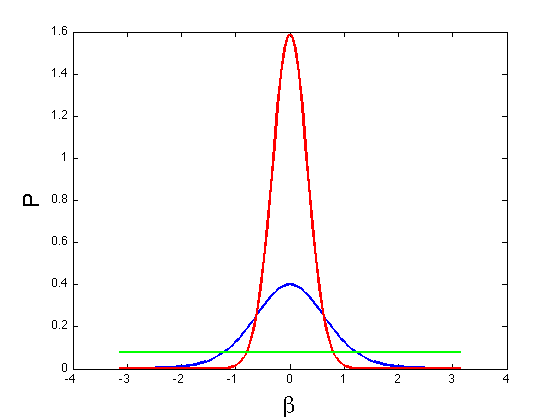,width=\linewidth,clip=}  \\[0.1cm] 
(a) \\
\epsfig{file=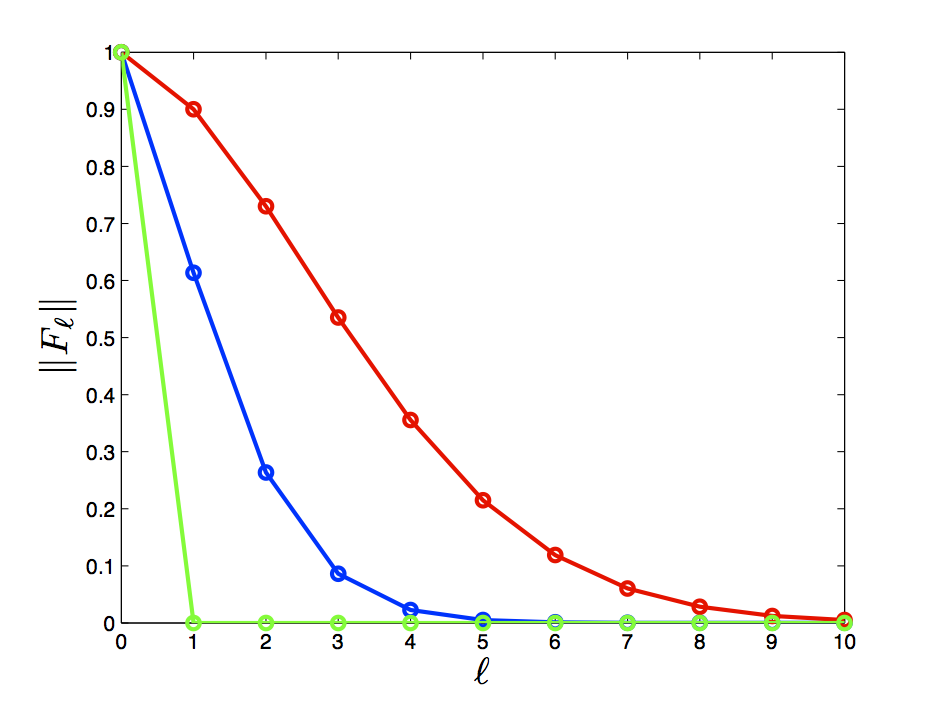,width=\linewidth,clip=} \\[0.1cm]
(b)
\end{tabular}
\caption{The Fisher-Von Mises distribution is shown in part (a) for $\kappa=0$ (green), $\kappa=2.5$ (blue), $\kappa=10$ (red).  The polar angle $\beta$ is shown from $-\pi$ to 
$\pi$ to center the distribution, though in practice only $[0,\pi]$ is used.  Part (b) shows the magnitude spectrum for the corresponding values of $\kappa$. }
\label{fig:fmvdist}
\end{figure}

A second example of spherical harmonic expansion is obtained for the Fisher Von-Mises distribution \cite{mardia}, which is defined with scale parameter $\kappa$ and mean $\mu \in S^{2}$ as
\beq
f(u; \mu, \kappa) = \frac{\kappa}{4\pi {\rm sinh} \kappa} e^{\kappa(\mu^{\top}u)}.
\label{eq:fishervm}
\eeq
Let $\mu = n$, the north pole. Then $n^{\top}u = \cos \beta$, and the distribution may be expanded in terms of the Legendre polynomials as follows \cite[eq 20]{schaeben}:
\beq
f(u;n,\kappa) = \frac{1}{4\pi} \sum_{\ell=0}^{\infty} (2\ell+1) \frac{I_{\ell+1/2}(\kappa)}{I_{1/2}(\kappa)} P_{\ell}(\cos \beta)
\label{eq:legnseries}
\eeq
Here, $I_\nu$ is the modified Bessel function of order $\nu$.  It is easily shown from this formula that $F_{\ell}^{m}=0$ unless $m=0$, and that
\beq
F_{\ell}^{0} = \| F_{\ell} \| =  \frac{I_{\ell+1/2}(\kappa)}{I_{1/2}(\kappa)}.
\label{eq:fvmcoeffs}
\eeq
Figure~\ref{fig:fmvdist}(a) and (b) show the distribution and its magnitude distribution for various $\kappa$,  

The magnitude-phase partition of spherical harmonic coefficient vectors, $F_\ell = \| F_{\ell} \| U_{\ell}$, invites us to consider which is more important to the structure -- magnitude or phase?  Considering that $F_{\ell}$ is $2\ell+1$ dimensional, but $\|F_\ell \|$ is a scalar, we should expect that the magnitude constrains only a small portion of the structure: phase should be much more important.  Suppose that $f$ is a real-valued function.  Then, due to conjugate symmetry, $F_{\ell}^{-n} = (-1)^n \left( F_{\ell}^{n} \right)^{*}$, so that each $F_{\ell}$ vector has only $2\ell+1$ real degrees of freedom.  Of these, one degree of freedom is constrained by the magnitude $\| F_\ell \|$, leaving $2\ell$ degrees of freedom in the phase vector $U_\ell = F_\ell / \| F_\ell \|$.  Consequently, for a function $f$ defined by its first $\ell\leq L$ spherical harmonic vectors, the magnitude spectrum $\|F_\ell\|$ for $\ell\leq L$ constrains only $L+1$ of the total
\beq
\sum_{\ell=0}^{L} 2\ell + 1 = (L+1)^{2}
\eeq
degrees of freedom. This means that the phase spectrum, $\{U_0, U_1, \ldots, U_{L}\}$ constrains the remaining  $L(L+1)$ degrees of freedom, or as a percentage
\beq
100 \frac{L}{L+1}.
\eeq
For example, for $L=10$, the phase spectrum constrains roughly $91\%$ of the signal.  Similar estimates can be derived for complex-valued signal. This is unlike the situation on $\field{R}$, where the phase and magnitude each constrain exactly half the degrees of freedom of the transform coefficients.  We see that on the sphere, phase is much more important to a signal than magnitude.

We now describe an important property of magnitude, that it only determines an axially-symmetric component of a function on $S^{2}$.  First, suppose that $f$ is axially-symmetric about the North pole $n=[0,0,1]$, i.e., $f(Ru)=f(u)$ for all rotations $R$ such that $Rn=n$.  Then, from (\ref{eq:rotS2}), we see that its spherical harmonic coefficient vectors satisfy $F_{\ell}=F_{\ell}D_{\ell}(R)$ for such rotations $R$.  From (\ref{eq:wignerdfunct}) we see that this is only possible if $F_{\ell}^{m} = 0$ for $m\neq 0$, i.e., $F_{\ell}$ has only one non-zero element, which located at the center of the $2\ell+1$ dimensional row vector.  Define the $2\ell+1$ row vector $Q_{\ell}$ with $1$ in its middle entry:
\beq
Q_{\ell} = [0,\ldots,0,1,0,\ldots,0]
\label{eq:pell}
\eeq
Therefore, functions axially-symmetric about the North pole satisfy $F_{\ell} = F_{\ell}^{0} Q_{\ell}$.  Consider now a spherical function $f$ about whose spherical harmonic coefficient vectors $\{ F_\ell \}_{\ell=0}^{\infty} $ we know only the magnitude spectrum $ \| F_{\ell} \| $.  Although there are infinitely many functions with the same magnitude spectrum, we see that one possibility is the axially-symmetric function whose spherical harmonic coefficient vectors satisfy $F_{\ell} = \| F_{\ell} \| Q_{\ell}$.   Since such functions are entirely determined by their magnitude spectra, we conclude that the magnitude component does not determine any more than an axially-symmetric version of the function.  

\subsection{The view from $SO(3)$}
\label{sec:3a}

The properties of spherical harmonics become clearer when $S^2$ is expressed as a homogeneous
space for the group $SO(3)$ of $3$-D rotations. That group consists of all $3\times 3$ orthogonal matrices with determinant $+1$.  
Every function $f$ whose domain is $S^2$ may be ``lifted'' to a corresponding function
$\widetilde{f}$ on $SO(3)$ by the mapping $\widetilde{f}(R) = f(Rn)$, where $n$ is the North pole as before. 

 Note that $\widetilde{f}$ is constant on rotations that fix $n$, and, therefore, 
\beq
\widetilde{f}(\alpha,\beta,\gamma)=\widetilde{f}(\alpha,\beta,0), \quad \forall \gamma.
\label{eq:tildprop}
\eeq
The elements $\{D_{\ell}^{mn}\}_{\ell,m,n}$ of the Wigner $D$-matrices mentioned above form an orthogonal basis for functions defined on $SO(3)$. We review some basic facts; see \cite{rockmoreffts}\cite{varshalovich} for details.  
The Fourier transform on $SO(3)$ consists of matrix-valued coefficients, one at each ``frequency'' $\ell$, where $\ell = 0,1,\ldots$. The
$\ell$-th Fourier coefficient is 
\beq
F(\ell) = \int_{SO(3)} f(R)D_{\ell}(R)^{\dagger} dR.
\label{eq:ftso3}
\eeq
Here the symbol $\dagger$ denotes matrix hermitian transpose and $dR =(8\pi^2)^{-1}\sin(\beta)d\alpha\,d\beta\,d\gamma$ is the normalized, rotation-invariant, measure on $SO(3)$. A $3$-D rotation of $f$, which equivalently translates the function on $SO(3)$, produces a corresponding transformation of the Fourier coefficients:  
\beq
g(R) = f(SR)\Leftrightarrow  G(\ell) =  F(\ell)D_\ell(S).
\label{eq:transso3}
\eeq
With (\ref{eq:ftso3}), the function $f$ may be expanded as follows
\beq
f(R) = \sum_{\ell=0}^{\infty} (2\ell+1){\rm Trace}\left[F(\ell) D_{\ell}(R)\right].
\label{eq:iftso3}
\eeq
The factor $(2\ell+1)$ corrects for the $L_{2}$ norm of the matrix elements:
\beq
\int_{SO(3)} \left| D_{\ell}^{mn}(R) \right|^2 dR = \frac{1}{2\ell+1}.
\eeq 


If $f = \widetilde{f}$, or in other words $f$ is a function lifted from $S^2$ as described above, then from
eqns.~(\ref{eq:wignerdfunct}) and (\ref{eq:tildprop}) we obtain that only the middle row of the $(2\ell +1)$-dimensional Fourier
coefficient matrix $\widetilde{F}(\ell)$ is non-zero.
\beq
m\neq 0 \Rightarrow \widetilde{ F}(\ell)^{mn} = 0.
\label{eq:middlerow}
\eeq
Moreover, from (\ref{eq:sphd}), we obtain (see \cite{kakaralacvpr}) that for lifted functions, the middle row of each of the $SO(3)$ Fourier coefficients contains, up to a scale factor, the Fourier coefficients
of the spherical harmonic expansion:
\beq
 \widetilde{ F}(\ell)^{0n}=\frac{1}{4 \pi c_{\ell}^{0}} { F}_{\ell}^{n}.  
\label{eq:middlerowc}
\eeq
It is easy to verify that, with (\ref{eq:middlerowc}) and (\ref{eq:sphd}), the expansion (\ref{eq:iftso3}) reduces to (\ref{eq:spharmexp}).  


 \section{Filtering}
 
Filtering, which is convolving a signal $f$ with a kernel $h$, is the most basic signal processing operation.  Due to the non-commutative nature of the rotation group $SO(3)$, filtering is not easy to define on the sphere $S^2$ in terms of an integration operation. In this section, we review several plausible definitions of filtering on the sphere, including a novel approach using lifting as defined in Section~\ref{sec:3a}.  Any useful definition of filtering, denoted $g=f\star h$, has the following properties.
\begin{description}
\item[F1] It maps spherical function $f$ to an output $g$ which is also a spherical function;
\item[F2] The kernel $h$ may be either axially-symmetric or directionally selective; 
\item[F3] It is associative, so that $(f\star h_1)\star h_2 = f\star(h_1\star_{o} h_2)$, where $\star_o$ is a suitable convolution on the domain of the filters $h_1$ and $h_2$;
\item[F4] It permits of FIR filtering, which is defined as a finite weighted sum of rotated copies of the function  
\beq
g(u) = \sum_{k=0}^{N-1} b_{k} f(R_{k}u)
\label{eq:fir}
\eeq 
The weights $\{b_{k}\}$ and the rotations $\{R_{k}\}$ define the operation of the filter. 
\end{description}
Note that simply rotating a function on $S^2$ is a filtering operation, and is a special case with $N=1$ of FIR filtering described in (\ref{eq:fir}).

We describe three approaches to filtering, and verify for each approach whether it has properties F1-F4.  The first approach, which we call the ``rotation'' approach, to filtering is to define convolution on $S^2$ analogously to that on the real line, so that the domain is $S^2$ and the translational variable, which in this case is rotation, is applied to the kernel
\beq
g(R) = (f\star h)(R) = \int_{S^{2}} f^{*}(u)h(Ru)du.
\label{eq:naive}
\eeq
This definition produces a function $g$ on $SO(3)$, not $S^{2}$, and thus does not have property F1 above. Although (\ref{eq:naive}) is not suitable, it is worth noting that it has a convenient expression in terms in the frequency domain.  The Fourier coefficients of $G$ on $SO(3)$ are outer-products of $F_\ell$ and $H_{\ell}$ as follows:
\beq
G(\ell) = \eta_{\ell} F_{\ell}^{\dagger} H_{\ell}
\label{eq:naivefreq}
\eeq 
with constants $\eta_{\ell}=1/({2\ell+1})$.  Note that both sides of the above are $2\ell+1$ dimensional matrices.  A proof using vector notation is provided in Appendix A;  a similar result using indices is shown by Kostelec \& Rockmore \cite{rockmoreffts}.  

The second definition of filtering, which may be called ``left convolution'',  is described and analyzed in Driscoll \& Healy \cite{Driscoll_Healy_1994}.  For $f$, $h$ in $L^{2}(S^2)$, and $n$ the North pole, the convolution is obtained by integrating over rotations:
 \beq
g(u) =  (f \star h) (u) = \int_{SO(3)} f(Rn) h(R^{-1} u) dR
 \label{eq:drishealsphere}
 \eeq
  It has been shown \cite{Driscoll_Healy_1994} that spherical harmonic coefficients of $g$ are obtained by 
 \beq
G_\ell =  \left( F\star H \right)_\ell = \frac{2\pi}{c_{\ell}^{0}}  F_{\ell} H_{\ell}^{0}.
\label{eq:drishealfilter}
 \eeq
 The coefficient $c_\ell^{0}$ is defined in (\ref{eq:cellzero}). Note that only the central coefficient of $H_{\ell}$ vector is involved in producing the output vector.  Therefore, the phase of $G_{\ell}$ is, up to sign, the phase of $F_{\ell}$, so that filtering defined by (\ref{eq:drishealsphere}) essentially preserves phase.  Since only the central coefficient is used, this definition only works for axially-symmetric kernels, and therefore does not have property F2.
 
The third definition of filtering, which is our proposal, is to convolve both functions $f$ and $h$ on $SO(3)$, and project the result back to the sphere. The basis of this definition is the convolution of two functions on $SO(3)$, which is defined as 
\beq
g(V) = \int_{SO(3)} h(R)f(R^{-1}V)dR.
\label{eq:convnative}
\eeq
Using (\ref{eq:ftso3}), it may be shown (see Appendix~B) that the Fourier transform of $g$ is related to those of $f$, $h$ as follows:
\beq
G(\ell) = F(\ell) H(\ell)
\label{eq:convso3nativefourier}
\eeq
We define the {\em projected convolution} on $SO(3)$ as follows.  Let $f$ be a function on $S^{2}$, and $h$ a kernel function on $SO(3)$, whose construction is discussed later.  Define 
\beq
g(u) = P_{S^{2}} \left[ \int_{SO(3)} h(R) \widetilde{f}(R^{-1}V) dR \right]
\label{eq:projconvonso3}
\eeq
Here the projection operator $P_{S^{2}}$ converts the three variable function on $SO(3)$ to a two-variable function on $S^{2}$ by integration.  Since $\alpha$, $\beta$ are sufficient to parameterize the sphere by (\ref{eq:unitvec}),  we integrate over the Euler angle $\gamma$ as follows:
\beq
P_{S^{2}}\left[g(R)\right] = \frac{1}{2\pi} \int_{0}^{2\pi} g(\alpha,\beta,\gamma) d\gamma.
\label{eq:s2}
\eeq 
The properties of (\ref{eq:projconvonso3}) are now established.   For each $\ell\geq 0$, let $G_{\ell}$ be  the $(2\ell+1)$-dimensional row vector of spherical harmonic coefficients for $g$, the filter output,  $F_{\ell}$ is the corresponding vector for $f$, and $H(\ell)$  is the square $(2\ell+1)$-dimensional Fourier coefficient matrix obtained from (\ref{eq:ftso3}) for the filter kernel $h$.

\begin{thm} With the notation as above, the projected convolution (\ref{eq:projconvonso3}) has the Fourier representation
\label{thm:filtering}
\beq
G_{\ell} = F_{\ell} H(\ell),
\label{eq:Fourierprojconv}
\eeq
Furthermore,  it satisfies all four properties F1-F4.   
\end{thm}
\begin{IEEEproof}
Using $Q_{\ell}$ as defined in (\ref{eq:pell}),  we have from (\ref{eq:middlerowc}) that 
\beq
\widetilde{F}(\ell) = \frac{1}{4\pi c_{\ell}^{0}} Q_{\ell}^{\top} F_{\ell}
\label{eq:widetfell}
\eeq
Let the integral in (\ref{eq:projconvonso3}) be denoted
\beq
g(V) = \int_{SO(3)} h(R) \widetilde{f}(R^{-1}V) dR.
\eeq
This is a function on $SO(3)$, and its Fourier coefficients, from (\ref{eq:convso3nativefourier}) and (\ref{eq:widetfell}) are
\beq
G(\ell) =  \frac{1}{4\pi c_{\ell}^{0}} Q_{\ell}^{\top} F_{\ell} H(\ell) 
\label{eq:gell}
\eeq
Consequently, using (\ref{eq:iftso3}), and rearranging terms, we see that (\ref{eq:projconvonso3}) becomes  
\beq
 \sum_{\ell=0}^{\infty} \frac{(2\ell + 1)}{4\pi c_{\ell}^{0}} {\rm Trace}\left[ Q_{\ell}^{\top} F_\ell H(\ell) \int_{0}^{2\pi} D_{\ell}(\alpha,\beta,\gamma) \frac{d\gamma}{2\pi} \right]
\eeq
Using (\ref{eq:wignerdfunct}) and (\ref{eq:sphd}), we obtain that
\beq
g(u) = \sum_{\ell=0}^{\infty} F_{\ell} H(\ell) Y_{\ell} (u)
\eeq
from which, by comparing with (\ref{eq:vecnotspharmexp}), the result (\ref{eq:Fourierprojconv}) follows.  

We now verify satisfaction of the properties.  F1 is obvious. To prove F2, note that the case of axially-symmetric filtering, which is provided by (\ref{eq:drishealfilter}), is obtained by setting $H(\ell) = \frac{2\pi}{c_\ell^0} H_{\ell}^{0} I_{\ell}$, where $I_{\ell}$ is the $2\ell+1$-dimensional identity matrix. The case of directionally-selective kernels is established by the FIR filtering property shown next. The property F3 follows from the transform (\ref{eq:Fourierprojconv}), since $(f\star h_1)\star h_2$ gives coefficients $F_{\ell} H_{1}(\ell) H_{2}(\ell)$, which by (\ref{eq:convso3nativefourier}) is the same result for transforming $f\star (h_1\star_o h_2)$.  To prove F4, note that the FIR filter is easily shown have the Fourier transform 
\beq
H(\ell) = \sum_{k=0}^{N-1} b_{k} D_{\ell}(R_{k})
\eeq 
Examples of directional FIR filters are given below.  
\end{IEEEproof}

The three types of filtering, including our proposal, are summarized, with a ``visualization''  in terms of vectors and matrices for the first coefficient $\ell=1$ in Table~\ref{tab:filtering}.

\begin{table}[h]
\begin{center}
    \begin{tabular}{| c| c| c| }
    \hline
   Equation & Domain of  $g$ & Illustration of $G$ for $\ell=1$ \\
   \hline
  (\ref{eq:naive}) & $SO(3)$ & $\left[\begin{array}{c}G\, G\, G \\ G\, G\, G \\ G \, G\, G \end{array} \right] \propto \left[\begin{array}{c}F\\ F\\ F \end{array}\right]\left[\begin{array}{c}H\,H\,H\end{array}\right]$ \\
  \hline
  & &  \\
  (\ref{eq:drishealsphere}) & $S^{2}$ & $\left[\begin{array}{c} G\,G\,G \end{array}\right] \propto \left[\begin{array}{c}F\, F\, F \end{array}\right]H$   \\    
  & &  \\
   \hline
   (\ref{eq:projconvonso3}) & $S^2$ & $\left[\begin{array}{c}G\, G\, G  \end{array} \right] = \left[\begin{array}{c} F\,F\,F\end{array}\right] \left[\begin{array}{c} H\, H\, H \\ H\, H\, H \\ H \, H\, H \end{array}\right]$ \\
   \hline
    \end{tabular}
    \caption{Summary of filtering methods  $g=f \star h$ on the sphere, with the proposed method on the bottom row. In the rightmost column, $G$, $F$, and $H$ represent generic (=indices suppressed) harmonic coefficients  of output, signal, and kernel, respectively.}
\label{tab:filtering}
\end{center}
\end{table}

Table~\ref{tab:filtering} helps to explain why (\ref{eq:projconvonso3}), our proposed method of filtering, is a ``phase-sensitive'' method.  The previously-known method of filtering on $S^2$, (\ref{eq:drishealsphere}), allows only axially-symmetric filters to be used.  Axially-symmetric filters, as shown in Section~\ref{sec:magphase}, are determined essentially by their magnitude spectrum and therefore have no phase component.  Furthermore, the effect of the previous method (\ref{eq:drishealsphere}) on the input function $f$ only changes its  magnitude components (up to sign) and not its phase.    In contrast, the proposed method (\ref{eq:projconvonso3}) allows freedom in designing the phase response of the filter, and consequently provides a wide range of filtering effects on both the magnitude and phase components of the input function. 

Here are some examples illustrating the projected convolution approach. The simplest filter is rotation: $g(u)=f(Ru)$ for some rotation $R$. The coefficients of this filter are unitary: $H(\ell) = D_{\ell}(R)$ for all $\ell$.  Therefore, the phase of the signal coefficients, $U_{\ell} = F_{\ell}/\|F_{\ell}\|$, is transformed linearly: $U_{\ell}\mapsto U_{\ell} D_{\ell}(R)$.   A second, relatively simple filter, is defined on $S^{2}$ and has an axially-symmetric kernel.  As discussed in the proof,  its $SO(3)$ Fourier coefficients are $H(\ell) = \frac{2\pi}{c_{\ell}^{0}} H_{\ell}^{0} I_{\ell}$, where $I_{\ell}$ is the $2\ell + 1$-dimensional identity matrix and $H_{\ell}^{0}$ are its central spherical harmonic coefficients.. 

A third example is a local average of  $f(u)$ and two neighbors at lines of higher and lower latitude, respectively.  This $3$-point filter may be described as 
\beq
g(u) = 0.5f(u) + 0.25f(R_{\beta_0}u) + 0.25f(R_{\beta_0}^\top u).
\label{eq:vertfilter}
\eeq
Here, $R_{\beta_0}$ is the rotation that moves the North pole $n$ to the axis $[\sin\beta_0,0,\cos\beta_0]$, which is aligned along the $\alpha=0$ meridian.  The corresponding filter has the Fourier transform
\beq
H({\ell}) = 0.5I_{\ell} + 0.25D_{\ell}(R_0) + 0.25D_{\ell}(R)^{\dagger}.
\label{eq:nsfilter}
\eeq
This is an example of a filter that is not axially symmetric, and which may be implemented using our projected convolution approach.  

A fourth and more detailed example, is obtained by applying the associativity property F3.   If we apply the FIR filter (\ref{eq:nsfilter}) to the output of an axially symmetric filter (such as the Fisher-Von Mises) then the result a weighted combination of symmetric filters arranged on the latitude axis.  Such a filter may be described as
\beq
g = 0.5f\star h + 0.25f\star h_{R} + 0.25f\star h_{R^{\top}}.
\label{eq:firweights}
\eeq
In the frequency domain, the filtering is described as $G_{\ell} = F_{\ell} H_{c}(\ell)$, where the cascaded filter has coefficients
\beq
H_{c}({\ell}) = \frac{2\pi}{c_{\ell}^{0}} H_{\ell}^{0} \left[0.5I_{\ell}+0.25D_{\ell}(R_{\beta_0})+0.25D_{\ell}(R_{\beta_0})^{\dagger}\right].
\label{eq:composite}
\eeq

We refer henceforth to the $H(\ell)$ matrices as the transfer function of the filter.  An important concept in filtering is the impulse response.  We obtain an analogous impulse response for the projected convolution using the inverse spherical harmonic expansion of (\ref{eq:Fourierprojconv}), where the input coefficient vectors $F_\ell$ are those for the spherical impulse as defined in (\ref{eq:impfourier}).   However, such an impulse response is incomplete, since it determines only a function on $S^2$ and not the kernel on $SO(3)$.  Hence, the transfer function is required as the complete description of the filter. 

\section{Experiments}

We implemented the projected convolution and validated its properties in three experiments.  In each experiment, we use the method of computing spherical harmonics described by Chung {\it et al} \cite{ChungDSED07}, and the software accompanying that paper. That method is called iterative residual fitting (IRF), which progressively computes at each frequency $\ell$, the spherical harmonic coefficients which best fit, in a least-squares sense, the residual obtained by subtracting from the data a weighted representation from the next lower frequency $\ell-1$. Intuitively, the progressive weighted coefficient fitting provides a ``windowing'' effect similar to that known to reduce Gibbs phenomenon for Fourier series on the real line.   We also used data and software accompanying the papers \cite{Driscoll_Healy_1994}\cite{rockmoreffts}\cite{simons06}.  MATLAB code for all of our experiments is available online \cite{onlinecode}.

\begin{figure*}
  \centering
  \subfloat[]{\includegraphics[width=0.3\linewidth]{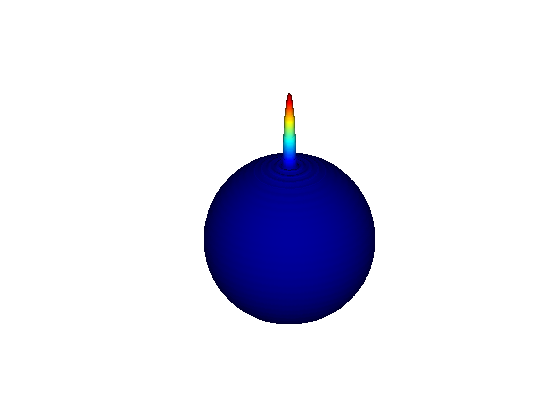}}
  ~ 
 \subfloat[]{\includegraphics[width=0.3\linewidth]{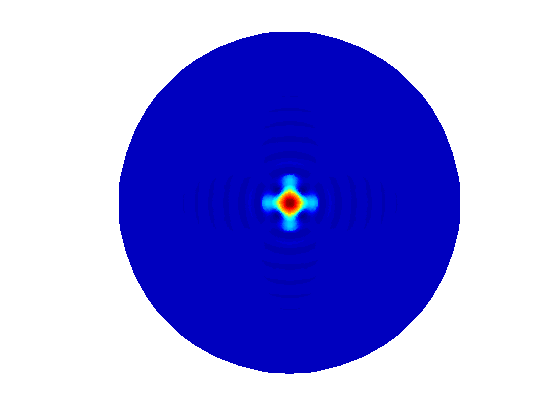}}
 ~ 
  \subfloat[]{\includegraphics[width=0.3\linewidth]{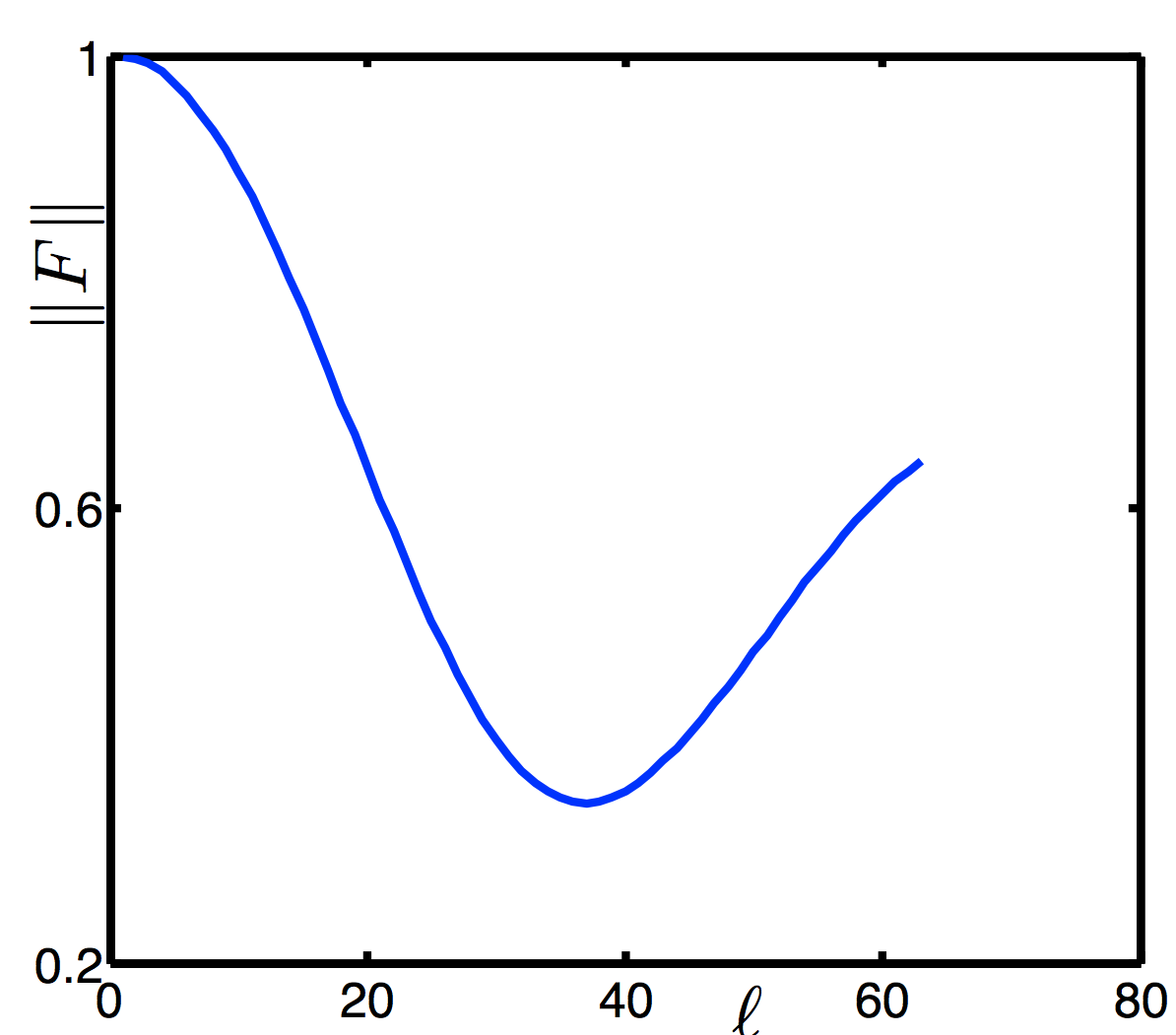}}
 \caption{(a) Shows the impulse function on a sphere reconstructed from its $64$-th order spherical harmonic coefficients; (b) shows the impulse response of the $5$-pt lowpass filter, with a central lobe and four sidelobes clearly visible; (c) shows the magnitude of the frequency response for the same filter, normalized by the magnitude spectrum of a single impulse to better illustrate its lowpass characteristic. }
  \label{fig:impulse}
\end{figure*}

\begin{figure*}
  \centering
  \subfloat[]{\includegraphics[width=0.3\linewidth]{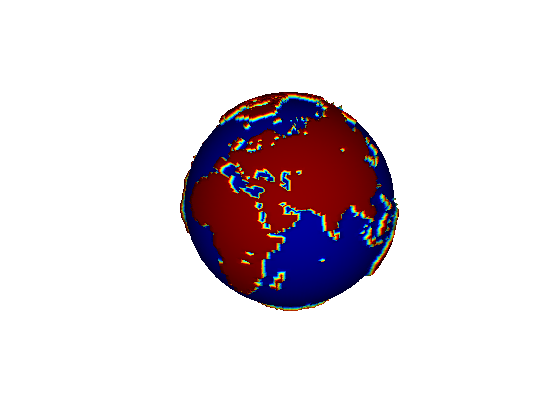}}
  ~ 
  \subfloat[]{\includegraphics[width=0.3\linewidth]{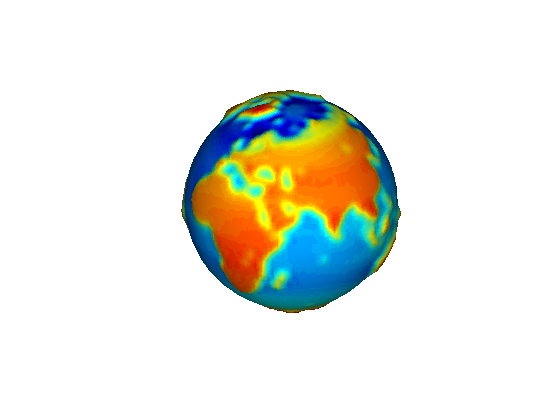}}
  ~ 
  \subfloat[]{\includegraphics[width=0.3\linewidth]{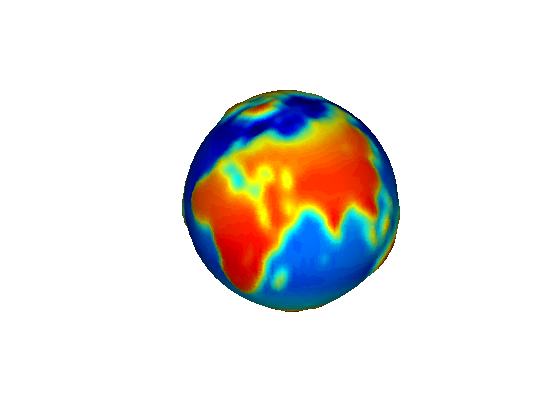}}
  \caption{Part (a) shows a binary world map, part (b) shows the reconstruction by spherical harmonics with bandwidth $63$, and part (c) shows the reconstruction filtered by the $5$-point lowpass filter in eq.~(\ref{eq:5ptfilter}), with the smoothing effect clearly visible.  }
  \label{fig:world}
\end{figure*}

In our first experiment, we construct a simple $5$-tap FIR lowpass filter.  Figure~\ref{fig:impulse} shows the reconstruction of an impulse defined in the frequency domain using (\ref{eq:impfourier}), using bandwidth $L=63$ as defined in (\ref{eq:finbandspharmexp}).  The reconstruction is obtained using IRF.  The filter is defined by the input-output equation
\begin{eqnarray}
g(u) &=& 0.5f(u) + 0.125 f(R_{1}u) + 0.125 f(R_{1}^{\top} u) + \nonumber \\ 
& & 0.125 f(R_{2} u) + 0.125 f(R_{2}^{\top} u)
\label{eq:5ptfilter}
\end{eqnarray}
The two rotation matrices, $R_1$ and $R_2$ are defined as follows: $R_1$ has Euler angles $\alpha=0$, $\beta=\pi/32$, and $\gamma=0$, while $R_2$ has Euler angles $\alpha=\pi/2$, $\beta=\pi/32$, $\gamma=-\pi/2$. The rotation $R_1$ shifts the North pole down by $\pi/32$ along the $0^{o}$ (Greenwich) meridian, and the second $R_2$ shifts by the same amount along the $90^{o}$E meridian.  This filter is easily implemented in the frequency domain using the coefficient matrices
\begin{eqnarray}
H(\ell) &=& 0.5I_{\ell} + 0.125\left[ D_{\ell}(R_1) +  D_{\ell}(R_1)^{\dagger}\right] + \nonumber \\ 
& & 0.125 \left[ D_{\ell}(R_2) +  D_{\ell}(R_2)^{\dagger}\right] 
\label{eq:5ptfilterfreq}
\end{eqnarray}
Figure \ref{fig:impulse} shows the impulse response of the filter, illustrating four smaller lobes surrounding a central main lobe.  The impulse response is computed using the inverse transform of the projected convolution, which is obtained in the frequency domain using (\ref{eq:Fourierprojconv}).  As a simple summary of the frequency response, we computed the norm $\| H_{\ell} \|$ at each frequency $\ell$, and plot the results in Figure~\ref{fig:impulse}.   The plot shows values after normalization by the magnitude of the impulse's coefficient vectors as shown in (\ref{eq:frac}).   It can be seen that the norm, which measures the ``size'' of the magnitude component, generally decreases with $\ell$, clearly the characteristic of a lowpass as expected from a local average.

\begin{figure*}
  \centering
  ~ 
  \subfloat[]{\includegraphics[width=0.3\linewidth]{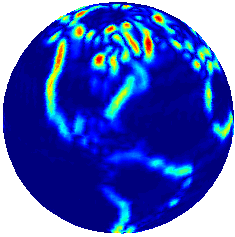}} \qquad
  ~ 
  \subfloat[]{\includegraphics[width=0.3\linewidth]{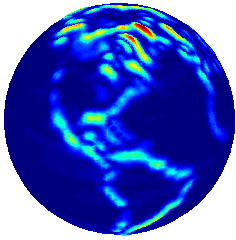}} \\ 
    \subfloat[]{\includegraphics[width=0.3\linewidth]{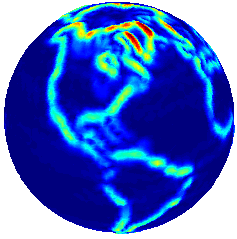}} \qquad
      \subfloat[]{\includegraphics[width=0.3\linewidth]{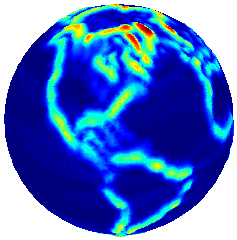}}
      \caption{Filtered versions of the world map.  Part (a) shows the result of filtering using the directional butterfly filter (\ref{eq:bflyplane}), which acts as an edge detector, and part (b) shows the result with the orthogonally-oriented filter, which detects edges in a perpendicular direction.  Part (c) shows the sum of absolute values of filter outputs in (a) and (b), illustrating that parts (a) and (b) form complementary aspects of an edge detector, and part (d) shows the output when the filter is dilated by a factor of two, showing a blurrier output as expected.  }
 \label{fig:worldbutterfly}
 \end{figure*}

In our second experiment, we filtered a binary world map shown in Figure~\ref{fig:world}.  In the same figure, we show a reconstruction of the world map using bandwidth $L=63$, as in (\ref{eq:finbandspharmexp}) but where the coefficients $F_\ell$ are obtained by IRF.  The value of IRF in this case is that it suppresses the Gibbs phenomenon ringing due to the binary edges of the world map.  We then applied the filter (\ref{eq:5ptfilterfreq}) to obtain the result in Figure~\ref{fig:world}(c).  The lowpass filtering of the $5$-tap filter is evident.  We then validated the theory further by using a directional filter.   The butterfly filter is defined by McEwan {\it et al} \cite{McEwenHML07} as follows: for $(x,y)$ in the plane $\field{R}^2$, let 
\beq
h(x,y) = x e^{-(x^2+y^2)/2\sigma}.
\label{eq:bflyplane}
\eeq
This filter may be mapped to the sphere using the co-latitude $\beta$ and longitude $\alpha$ as 
\beq
h(\beta,\alpha) = \left[ \tan(\beta/2) \cos \alpha\right] e^{-\tan^2(\beta/2)/2\sigma}
\label{eq:betaalpha}
\eeq
We generate an FIR filter from this prototype by using a set of $N=144$ evenly-spaced samples on an angular grid in $\beta$, and $\alpha$.  To determine an FIR filter we need a third Euler angle $\gamma$, in order to set the rotation.  We set  $\gamma = -\alpha$ to counter the rotation around $z$ by $\alpha$, and let $R_k$, for $k=1$ to $N$, be the rotation defined by Euler angles $\alpha_k$, $\beta_k$, $\gamma_k=-\alpha_k$.  Hence, the filter is specified on the sphere by the input-output equation
\beq
g(u)  = \sum_{k=0}^{N-1} h(\beta_k,\alpha_k) f(R_k u).
\eeq
In the frequency domain, the filter is specified at each frequency $\ell$ by the transfer function
\beq
H(\ell) = \sum_{k=0}^{N-1} h(\beta_k,\alpha_k) D_{\ell}(R_k).
\label{eq:transferell}
\eeq
We use the transfer function in our implementation, constructing $H(\ell)$ for $\ell=0,\ldots,63$.   This filter is easily dilated by scaling each angle with a constant $\lambda$ as follows: $\beta_k \mapsto \lambda \beta_k$, $\alpha_k  \mapsto \lambda \alpha_k$.  This relative ease of dilation suggests multi-resolution implementations, though we do not have space to pursue that here.  Dilation may also be accomplished through stereographic projection \cite{Antoine97waveletson} \cite{Antoine1999}\cite{Wiauxeuclid} or through harmonic scaling 
\cite{wiaux2007sdw}.   Figure~\ref{fig:worldbutterfly} shows the results of filtering the world map with the FIR approximation to the butterfly, and the $90$-degree rotation of it using $y$ in place of $x$ in (\ref{eq:bflyplane}).  Also shown are the results using a dilated filter with $\lambda = 2$.  

 In our third experiment, we extended the FIR filtering method proposed in this paper beyond the sphere $S^2$ to a larger set of two-dimensional manifolds.  Let $f:S^2\rightarrow \field{R}^{3}$ be mapping with components denoted $f_x$, $f_y$, and $f_z$.  The components define the $x$, $y$, $z$ coordinates of points on the surface of the manifold.  Each of the three coordinate functions may be expanded using spherical harmonics.  For each $\ell$, let ${\cal F}_{\ell}$ denote the $3\times (2\ell + 1)$-dimensional matrix whose $3$ rows contain the spherical harmonic coefficients of the corresponding coordinate function.  Then the manifold equivalent of (\ref{eq:vecnotspharmexp}) is
\beq
 f(u) = \sum_{\ell=0}^{\infty} {\cal F}_{\ell} Y_{\ell}(u).
\eeq
This expansion is referred to in the computer vision literature as SPHARM, the capital letters signifying that it is a 3-D version of spherical harmonics \cite{BrechbuhlerGK95}.  FIR filtering may be performed on SPHARM coefficients with a simple extension of (\ref{eq:Fourierprojconv}), resulting in the equation
\beq
{\cal G}_{\ell} = {\cal F}_{\ell} H(\ell)
\label{eq:vecFourierprojconv}
\eeq

\begin{figure*}
  \centering
  \subfloat[]{\includegraphics[width=0.32\linewidth]{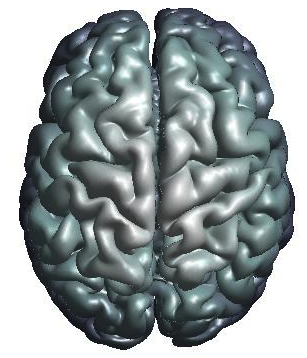}}
    \subfloat[]{\includegraphics[width=0.3\linewidth]{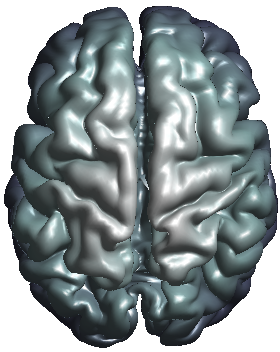}}
      \subfloat[]{\includegraphics[width=0.3\linewidth]{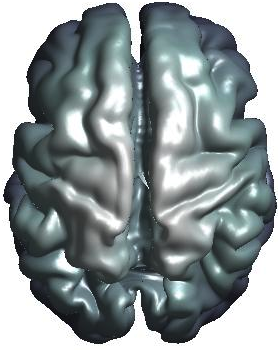}}

      \caption{Filtered versions of cortical surface. The manifold in part (a) is filtered in the SPHARM domain using the lowpass filter in eq. (\ref{eq:5ptfilter}), giving the result shown in (b), and with a cascade of the same filter applied twice, giving a twice smoothed result as shown in part (c).   }
 \label{fig:brain}
 \end{figure*}

Note that in the above equation, ${\cal G}$ and ${\cal F}$ are $3\times (2\ell+1)$, while the transfer function matrix $H(\ell)$ is $2\ell + 1$-dimensional.  Consequently, every filter developed for spherical harmonics may also be applied to SPHARM using exactly the same transfer function.  Figure~\ref{fig:brain} shows the result of filtering the expansion of a cortical surface using the lowpass filter of (\ref{eq:5ptfilterfreq}), as well as the result of filtering with the product $H(\ell)H(\ell)$ that is equivalent to cascaded filtering.  

 In our experiments, computation time is dominated by the IRF, which takes approximately $16$ seconds on a desktop computer with quad core $2.6$GHz CPU for a $L=64$ degree expansion of the world map data. Other steps, such the filtering operation (\ref{eq:Fourierprojconv}), are negligible in comparison.  The IRF may be replaced by the fast spherical harmonic transform  \cite{McEwenW11}  in cases of data where the Gibbs phenomenon caused by discontinuities is not a concern.  
 
\section{Summary and Conclusions}

One of the key developments in the research leading to this paper is relatively simple: writing the spherical harmonic expansion in vector form (\ref{eq:vecnotspharmexp}).  Once the expansion is seen in terms of vectors, the definition of phase, and its properties  in relation to magnitude become tractable.   It also becomes clear that to obtain a new set of vectors in a linear fashion, it is necessarily to multiply the existing set with compatible matrices.  The vector-matrix form of the filtering described in Theorem~\ref{thm:filtering} shows how such filtering is constructed.  We have shown the results of experiments that validate the theory using both spherical and manifold data.  

One conclusion we reach in this work is that phase is more important than magnitude for spherical harmonics, and that a study of phase leads to better understanding of filtering.   A second conclusion is that it is important to use vector-matrix notation, reducing the clutter of indices, to understand the operations involved in filtering.  Finally, we conclude that FIR filtering on the sphere may be studied using the same formal methods of magnitude-phase analysis and transfer function derivation that serve so well in the familiar real-line case.  This last conclusion suggests many opportunities for future work, to translate FIR filtering theory, including filter design and filter banks, to the spherical case using methods similar to those employed in this paper.  

\section*{Appendix A}
\label{sec:appa}
We prove (\ref{eq:naivefreq}).  Substituting from (\ref{eq:vecnotspharmexp}) into (\ref{eq:naive}), and using the rotation property (\ref{eq:rotS2}), we obtain
\beq
g(R) = \sum_{\ell} \sum_{k} \int_{S^{2}} F_{\ell}^{*} Y_{\ell}(u)^{*} H_{k} D_{k} (R) Y_{k}(u) du.
\eeq
Since the integrand is a scalar, we apply the trace and use its linearity and cyclic permutation invariance to write 
\beq
g(R) = \sum_{\ell} \sum_{k} {\rm Trace}\left[ \int_{S^{2}} Y_{\ell}(u)^{*} A_k(R) Y_{k}(u) du F_{\ell}^{*} \right] 
\eeq
Here $A_k(R) = H_k D_k(R)$ is a row-vector.  Due to orthonormality (\ref{eq:normal}), the integral vanishes if $\ell \neq k$, and for $\ell=k$ becomes the {\em column} vector $A_{\ell}^\top$.  Consequently, we have 
\beq
g(R) = \sum_{\ell} {\rm Trace}\left[ D_k(R)^{\top} H_{k}^\top F_{\ell}^{*} \right] 
\eeq
The result now follows by noting that ${\rm Trace}(X) = {\rm Trace}(X^{\top})$, and comparing with eq. (\ref{eq:iftso3}).  

\section*{Appendix B}
\label{sec:appb}
We prove (\ref{eq:convso3nativefourier}).  From ({\ref{eq:ftso3}) we have 
\[
G(\ell) = \int_{SO(3)} g(V) D_{\ell}(V)^{\dagger} dV.
\]
Substituting from (\ref{eq:convnative}), we get (all integrals are over SO(3)) that 
\[
G(\ell) = \int\left[ \int h(R)  f(R^{-1}V)dR\right] D_{\ell}(V)^{\dagger} dV
\]
Reversing the order of integrals, we obtain
\[
G(\ell) = \int h(R) \left[ \int f(R^{-1}V) D_\ell(V)^{\dagger} dV \right] dR.
\]
Let $U=R^{-1}V$, so that $RU=V$.  Note that $dV$ is a Haar measure on a compact group, so the integral is invariant to left shift (see for example \cite{HewittRossII}).  Consequently.
\[
G(\ell) = \int h(R) \left[ \int f(U) D_\ell(RU)^{\dagger} dU \right] dR.
\]
Now $D_\ell(RU)^{\dagger} = D_\ell(U)^{\dagger} D_\ell(R)^{\dagger}$, so that we obtain
\[
G(\ell) = \int h(R) \left [ \int f(U) D_\ell(U)^{\dagger} dU \right] D_\ell(R)^{\dagger} dR.
\]
The term in brackets is $F(\ell)$ by (\ref{eq:ftso3}), and hence (\ref{eq:convso3nativefourier}) follows.
\[
G(\ell) = \int h(R) F(\ell) D_{\ell}(R)^{\dagger}dR = F(\ell) H(\ell).
\]

\section*{Acknowledgements}
We thank W. Li, P. Kaliamoorthi, V. Premachandran, as well as the anonymous reviewers, for helpful comments.  

\bibliographystyle{ieee}
\bibliography{sphericalfiltering}

\end{document}